\documentclass[11pt]{article}
\usepackage{amsmath,amssymb,amsthm,amsfonts}
\usepackage{cite}
\usepackage[hidelinks]{hyperref}
\allowdisplaybreaks[4]
\newtheorem{thm}{Theorem}
\newtheorem{lem}{Lemma}[section]

\newtheorem{prop}[lem]{Proposition}
\theoremstyle{remark}

\numberwithin{equation}{section}
\makeatletter
\renewenvironment{proof}[1][\proofname]{\par
  \pushQED{\qed}%
  \normalfont\topsep6\p@\@plus6\p@\relax
  \trivlist
  \item[\hskip\labelsep\normalfont\bfseries #1\@addpunct{.}]\ignorespaces
}{%
  \popQED\endtrivlist\@endpefalse
}
\makeatother
\newcommand{\R}{\mathbb{R}}
\newcommand{\norm}[1]{\left\lVert #1\right\rVert}
\DeclareMathOperator{\dist}{dist}
\DeclareMathOperator{\supp}{supp}
\DeclareMathOperator{\diver}{div}
\begin{document}
\title{Uniqueness and Nonuniqueness of Positive Bounded Solutions for
Semilinear Dirichlet Problems in Epigraphs}
\author{ Wei Ji\thanks{School of Mathematics and Statistics, and  Key Laboratory of Nonlinear Analysis $\&$ 
Applications (Ministry of Education), Central China Normal University, P.O. Box 71010, Wuhan 430079, P. R. China. 
Email: \texttt{jiwei2020@amss.ac.cn}. }
}
\date{ }
\maketitle

\begin{abstract}
We study the uniqueness and nonuniqueness of positive bounded solutions
for $-\Delta u=f(u)$ in uniformly $C^{2,\gamma}$ epigraphs
$\Omega\subset\R^d$, $d\geq2$, with $u=0$ on $\partial\Omega$.
We prove the uniqueness of positive bounded solutions for any positive reaction in dimension two.
In each dimension $d\geq3$, we construct positive reactions
and epigraphs
such that this equation
admit multiple positive bounded solutions.
Our results also address Open Question~2 of
[H.~Berestycki and C.~Graham, ARMA, 2026]
on uniqueness in epigraphs.
\end{abstract}

\medskip

\section{Introduction}
\label{sec-introduction}

In this paper, 
we study the existence and uniqueness of positive bounded solutions for the
following Dirichlet problem
\begin{equation}\label{mainequation}
 \begin{cases}
 -\Delta u=f(u)\,\ \ \text{in}\ \ \Omega,\\
 u=0\,\ \ \text{on}\ \ \partial\Omega,
 \end{cases}
\end{equation}
where $d\geq2$ and $\Omega\subset\R^d$ is an epigraph of the form
\begin{equation}\label{eq-epigraph-domain}
 \Omega=\{x=(x',y)\in\R^{d-1}\times\R:y>\varphi(x')\},\ \ 
 \varphi\in C(\R^{d-1}).
\end{equation}
Following Berestycki and Graham \cite{BG2025,BG2026},
we assume that $\Omega$ is uniformly $C^{2,\gamma}$
for some $\gamma\in(0,1]$, in the sense that after
a translation and a rotation at each boundary point,
$\partial\Omega$ is locally the graph of a function
on a ball of fixed radius, with a $C^{2,\gamma}$ norm
bounded independently of the boundary point{,
see} \cite[Definition~A.1]{BG2025}.
We assume that $f$ satisfies
\begin{equation}\label{eq-reaction}
 \begin{gathered}
 f\in C^{1,\gamma}_{\mathrm{loc}}([0,+\infty)),
 \ \ f(0)=f(1)=0,\ \ f'(0)>0,\ \ f'(1)<0,\\
 f(s)>0\ \ \text{for any}\ \ 0<s<1,
 \ \ f(s)<0\ \ \text{for any}\ \ s>1.
 \end{gathered}
\end{equation}
As in \cite{BG2026},
we call $f$ satisfying \eqref{eq-reaction} a positive reaction.

{
Problem \eqref{mainequation} describes the steady states
of the following reaction--diffusion equation with zero
Dirichlet boundary conditions
\begin{equation}\label{eq-parabolic}
 \begin{cases}
  \partial_t U=\Delta U+f(U)
  \,\ \ \text{in}\ \ (0,+\infty)\times\Omega,\\
  U=0
  \,\ \ \text{on}\ \ (0,+\infty)\times\partial\Omega.
 \end{cases}
\end{equation}
Such equations are widely used to model population growth
and spatial propagation
\cite{Fisher1937,KPP1991,AronsonWeinberger1978}.
In population models, $U=U(t,x)$ denotes the population
density at time $t$ and position $x$, the diffusion term
$\Delta U$ describes spatial dispersal, and $f(U)$
describes local population growth
\cite{CantrellCosner2003}.
The zero Dirichlet boundary condition models absorption
at the habitat boundary \cite{BG2022,GH2012}.
A stationary solution of \eqref{eq-parabolic} has the form $U(t,x)=u(x)$
and satisfies $\partial_tU=0$, so
\eqref{eq-parabolic} reduces to \eqref{mainequation}.

On bounded smooth domains, Brezis and Oswald \cite{BO1986}
proved the uniqueness of positive solutions for
\eqref{mainequation} when $s\mapsto f(s)/s$ is strictly
decreasing on $(0,+\infty)$.
For positive reactions, the additional requirement that
$s\mapsto f(s)/s$ be strictly decreasing on $(0,1]$ is
called the strong-KPP condition.
Under the strong-KPP condition, Berestycki and Graham \cite{BG2025}
proved the uniqueness of positive bounded solutions of
\eqref{mainequation} on uniformly $C^{2,\gamma}$ domains
$\Omega$, possibly unbounded, when $f'(0)$ lies outside
the closure of the principal limit spectrum of the
Dirichlet Laplacian on $\Omega$.

For $f$ satisfying \eqref{eq-reaction}, the geometry
of $\Omega$ also plays a role in the uniqueness
of positive bounded solutions for \eqref{mainequation}.
When $f$ satisfies \eqref{eq-reaction},
$\R^d\setminus\Omega$ is strongly star-shaped, and
$(\R^d\setminus\Omega)\setminus B_R(0)
=C\setminus B_R(0)$ for some $R>0$ and some convex set
$C\subset\R^d$, Berestycki and Graham
\cite{BG2026} proved that \eqref{mainequation}
has a unique positive bounded solution.
For $f$ satisfying \eqref{eq-reaction} and uniformly
$C^{2,\gamma}$ domains with a star-shaped complement,
we proved in \cite{JW2026} the uniqueness and nondegeneracy
of positive bounded solutions for \eqref{mainequation}.
When $f$ satisfies \eqref{eq-reaction} and $\Omega$ is
a uniformly Lipschitz epigraph, Berestycki, Caffarelli and Nirenberg
\cite{BCN1997} used the sliding method
\cite{BN1991} to prove the uniqueness
of positive bounded solutions for \eqref{mainequation}.
For bistable reactions with $\int_0^1f(s)\,ds>0$
and uniformly Lipschitz epigraphs bounded from below,
Berestycki, Graham and Wei \cite{BGW2025}
proved the uniqueness of positive bounded solutions
for \eqref{mainequation}.

For bistable reactions with $\int_0^1f(s)\,ds>0$,
such as 
\[
f(s)=s(s-\theta)(1-s), \ \ 0<\theta<1/2,
\]
Berestycki, Graham and Wei \cite{BGW2025}
proved the uniqueness of positive bounded solutions
for \eqref{mainequation} in uniformly Lipschitz epigraphs
bounded from below.
Under \eqref{eq-reaction}, Berestycki and Graham
\cite{BG2026} proved the uniqueness of positive bounded
solutions for \eqref{mainequation} in asymptotically
uniformly Lipschitz epigraphs.
This result requires an additional geometric
condition at infinity and does not cover all uniformly
$C^{2,\gamma}$ epigraphs.
Berestycki and Graham
\cite[Section~8.1, Open Question~2]{BG2026}
therefore posed the following open question:

\emph{If $f$ is positive, does
\eqref{mainequation} have exactly one positive bounded
solution on all epigraphs?}

We prove that the answer is affirmative in dimension $d=2$
and give counterexamples in each dimension $d\geq3$
for a class of positive reactions.
Our first main result is the following theorem,
which also allows $f'(0)=0$.
}

\begin{thm}\label{thm-main}
Suppose that $\gamma\in(0,1]$ and $f$ satisfies
\eqref{eq-reaction} with $f'(0)\geq0$
in place of $f'(0)>0$.
Then:
\begin{enumerate}
 \item[(i)]
 For $d\geq2$, \eqref{mainequation}
 admits a positive bounded solution in any uniformly
 $C^{2,\gamma}$ epigraph $\Omega\subset\R^d$.
 If $d=2$ and there exist $c,s_0>0$ such that
 \begin{equation}\label{eq-planar-quadratic}
  f(s)\geq cs^2
  \ \ \text{for any}\ \ 0<s<s_0,
 \end{equation}
then the positive bounded solution of \eqref{mainequation} is unique.

 \item[(ii)]
 If the function $q(s):=f(s)/s$ is concave on $(0,1)$
 and $q(s)>f'(0)$ for some $s\in(0,1)$,
 then for $d\geq3$, there exists
 a uniformly $C^{2,\gamma}$ epigraph $\Omega\subset\R^d$
 such that \eqref{mainequation} admits two positive
 bounded solutions $u$ and $U$ satisfying
 $0<u<U<1$ in $\Omega$.
\end{enumerate}
\end{thm}

{
We remark that \eqref{eq-planar-quadratic} holds
for any $f$ satisfying \eqref{eq-reaction}, and hence
Theorem~{\ref{thm-main}}(i) yields the uniqueness of positive
bounded solutions for \eqref{mainequation} in any uniformly
$C^{2,\gamma}$ epigraph $\Omega\subset\R^2$
under weaker assumptions on $f$
than those in \eqref{eq-reaction}.
In particular, when $f'(0)=0$, the condition
\eqref{eq-planar-quadratic} is satisfied if
$f(s)\sim As^\beta$ as $s\searrow0$
with $A>0$ and $1<\beta\leq2$.
This also provides a partial answer to
\cite[Section~8.4, Open Question~6]{BG2026}
for planar epigraphs.
Theorem~{\ref{thm-main}}(ii) shows that even when
$f$ satisfies \eqref{eq-reaction} and the additional
assumptions on $q(s)=f(s)/s$, the uniqueness of positive
bounded solutions for \eqref{mainequation} may fail
in uniformly $C^{2,\gamma}$ epigraphs
$\Omega\subset\R^d$ with $d\geq3$.
An calssical example is
the cubic reactions
\[
 f(s)=s(s+a)(1-s),\ \ 0<a<1.
\]

}

{
The construction for Theorem~\ref{thm-main}(ii)
also gives the following result.

\begin{thm}\label{thm-extension}
Suppose $s_0>0$ and $g\in C^{1,\gamma}([0,s_0])$ satisfy
$g(0)=0$ and $g(s)>0$ for any sufficiently small $s>0$.
Then there exist $\delta\in(0,s_0)$ and a reaction $f$
satisfying \eqref{eq-reaction} with $f'(0)=g'(0)\geq0$
in place of $f'(0)>0$ such that $f=g$ on $[0,\delta]$.
For each integer $d\geq3$, there exists a uniformly
$C^{2,\gamma}$ epigraph $\Omega\subset\R^d$
 on which
\eqref{mainequation} admits two positive bounded solutions
$0<u<U<1$.
\end{thm}
}

{
Our uniqueness arguments use the smallest and largest
positive bounded solutions $\underline u$ and $\overline u$
of \eqref{mainequation}, which satisfy
$\underline u\leq u\leq\overline u$ in $\Omega$ for any
positive bounded solution $u$ of \eqref{mainequation}.
In dimension two, we combine the strict vertical
monotonicity of $\underline u$ and $\overline u$
with the Hopf lemma and a compactness argument for
translated solutions to prove that
$\underline u\equiv\overline u$ in $\Omega$.

We use the method of upper and lower solutions and monotone
iteration to prove the nonuniqueness results.
For Theorem~\ref{thm-main}(ii), we use the solutions
on a strip obtained in \cite[Theorem~1.1]{DRR2025}
to construct two positive bounded solutions
of \eqref{mainequation}: along the same horizontal line,
one tends to zero at infinity, whereas the other remains
bounded away from zero.
}

{
The paper is organized as follows.
In Section~\ref{sec-extremal}, we prove
Theorem~\ref{thm-main}(i).
In Section~\ref{sec-high}, we complete the proofs of
Theorem~\ref{thm-main}(ii) and Theorem~\ref{thm-extension}.
}

\section{Uniqueness for \texorpdfstring{$d=2$}{d=2}}\label{sec-extremal}

{
The purpose of this section is to prove
Theorem~\ref{thm-main}(i) using the
minimal and maximal solutions.
}

Let $B_R(x)$ be the open ball of radius $R$ centered
at $x$ and set $B_R:=B_R(0)$.
We first establish bounds for positive bounded solutions
of \eqref{mainequation} on a uniformly $C^{2,\gamma}$
domain $G\subset\R^d$.
We denote by $H^1_{\mathrm{loc}}(\overline G)$ the space
of functions $v$ satisfying $v\in H^1(G\cap B_R)$
for any $R>0$.
Boundary values are understood in the trace sense,
and any $v\in H^1_{\mathrm{loc}}(\overline G)$ with zero
trace on $\partial G$ extends by zero to a function
in $H^1_{\mathrm{loc}}(\R^d)$.

\begin{lem}\label{lem-domain-bounds}
Let $G\subsetneq\R^d$ be a connected uniformly $C^{2,\gamma}$ domain. Suppose
that $f$ satisfies \eqref{eq-reaction}, or that $d=2$, $G=\Omega$ is an
epigraph as in \eqref{eq-epigraph-domain} and $f$ satisfies
\eqref{eq-reaction} with $f'(0)>0$ replaced by \eqref{eq-planar-quadratic}.
For each $\theta\in(0,1)$, each positive bounded solution $u$ of
\eqref{mainequation} on $G$ satisfies, for some $R_\theta>0$,
\begin{equation}\label{eq-domain-bounds}
 0<u<1\,\ \ \text{in}\ \ G,\ \  u\geq\theta\,\ \ \text{on}\ \ \{x\in
 G:\dist(x,G^c)\geq R_\theta\}.
\end{equation}
\end{lem}

\begin{proof}
It follows from \cite[Lemma~2.2]{JW2026} that
\eqref{eq-domain-bounds} holds when $f'(0)>0$.
Suppose $d=2$, $G=\Omega$ and $f'(0)=0$.
We first prove that $0<u<1$ in $\Omega$.
Define $z_+:=\max\{z,0\}$ for any $z\in\R$.
Since $u=0$ on $\partial\Omega$ and $f(s)<0$
for any $s>1$, the zero extension of $(u-1)_+$
is bounded and subharmonic in $\R^2$.
Applying the planar Liouville theorem and the strong
maximum principle, we then obtain that $0<u<1$ in $\Omega$.

We next determine the possible limits of $u(x)$
as $\dist(x,\Omega^c)\to+\infty$.
For any $R>0$, we obtain from \eqref{eq-epigraph-domain} that
\begin{equation*}
 \begin{aligned}
 \Omega_R
 &:=\{x\in\Omega:\dist(x,\Omega^c)>R\}\\
 &=\left\{(x',y):
 y>\max_{|t|\leq R}
 \left(\varphi(x'+t)+\sqrt{R^2-t^2}\right)\right\}.
 \end{aligned}
\end{equation*}
Suppose $\{x_n\}\subset\Omega$ satisfies
$\dist(x_n,\Omega^c)\to+\infty$ as $n\to\infty$.
Applying interior elliptic estimates, we then deduce
from \eqref{mainequation} that, up to a subsequence,
\begin{equation*}
 \begin{gathered}
 u(x_n+\cdot)\to U
 \,\ \ \text{in}\ \ C^2_{\mathrm{loc}}(\R^2)
 \ \ \text{as}\ \ n\to\infty,\\
 0\leq U\leq1,\ \ -\Delta U=f(U)\geq0
 \,\ \ \text{in}\ \ \R^2.
 \end{gathered}
\end{equation*}
Since $U$ is bounded and superharmonic in $\R^2$,
the planar Liouville theorem implies that $U$ is constant.
Hence $U\equiv0$ or $U\equiv1$.
It follows that for each $\varepsilon\in(0,1/2)$,
\[
 u(\Omega_R)\cap[\varepsilon,1-\varepsilon]
 =\varnothing
\]
for any sufficiently large $R>0$.
Since $\Omega_R$ is connected and decreases with $R$,
we conclude that
\begin{equation}\label{eq-deep-alts}
  \text{either}\ \
 \sup_{\Omega_R}u\to0
 \ \ \text{or}\ \ 
 \inf_{\Omega_R}u\to1
 \ \ \text{as}\ \ R\to+\infty.
\end{equation}

{\setlength{\emergencystretch}{1em}
We now exclude the first alternative in \eqref{eq-deep-alts}.
On the contrary, suppose that
$\sup_{\Omega_R}u\to0$ as $R\to+\infty$.
If $\varphi$ is bounded above on a horizontal half-line,
then $\Omega_R$ contains a translated upper quadrant
for any $R>0$.
Choose $R>0$ sufficiently large such that $u<s_0$
in $\Omega_R$, where $s_0$ is as in
\eqref{eq-planar-quadratic}, and denote by $Q$
a translated upper quadrant contained in $\Omega_R$.
We deduce from \eqref{mainequation}
and \eqref{eq-planar-quadratic} that
\begin{equation*}
 -\Delta u=f(u)\geq cu^2
 \,\ \ \text{in}\ \ Q,
\end{equation*}
}
which contradicts \cite[Theorem~5.1]{AS2011}.
We thus obtain that
\[
 \limsup_{x'\to\pm\infty}\varphi(x')=+\infty.
\]
For any $a\in\R$, each connected component $D$
of $\Omega\cap\{y<a\}$ is bounded.
Define for any $\lambda<a$,
\[
 D_\lambda:=D\cap\{y<\lambda\},\ \ 
 v_\lambda(x',y):=u(x',2\lambda-y)-u(x',y)
 \,\ \ \text{in}\ \ D_\lambda.
\]
It follows from \eqref{mainequation} that
\begin{equation*}
 \begin{gathered}
 -\Delta v_\lambda
 =v_\lambda\int_0^1f'(u+tv_\lambda)\,dt
 \,\ \ \text{in}\ \ D_\lambda,\\
 v_\lambda>0
 \,\ \ \text{on}\ \ 
 \partial D_\lambda\cap\{y<\lambda\},\ \ 
 v_\lambda=0
 \,\ \ \text{on}\ \ 
 \partial D_\lambda\cap\{y=\lambda\}.
 \end{gathered}
\end{equation*}
The small-domain maximum principle
\cite[Proposition~1.1]{BN1991} implies that $v_\lambda\geq0$
in $D_\lambda$ when $|D_\lambda|$ is sufficiently small.
Define
\[
 \lambda_*:=\sup\left\{
 \lambda<a:
 v_\mu\geq0\text{ in }D_\mu
 \text{ for any }\mu\leq\lambda
 \right\}.
\]
Suppose $\lambda_*<a$.
The continuity of $u$ yields that
$v_{\lambda_*}\geq0$ in $D_{\lambda_*}$.
Since $v_{\lambda_*}>0$ on
$\partial D_{\lambda_*}\cap\{y<\lambda_*\}$,
the strong maximum principle implies that
$v_{\lambda_*}>0$ in each connected component
of $D_{\lambda_*}$.
Since $D$ is bounded, we obtain from the dominated
convergence theorem \cite[Proposition~1.1]{BN1991} that
\begin{equation*}
 |\{x\in D_\lambda:v_\lambda(x)<0\}|\to0
 \ \ \text{as}\ \ \lambda\searrow\lambda_*.
\end{equation*}
Applying the small-domain maximum principle on
$\{x\in D_\lambda:v_\lambda(x)<0\}$,
we then obtain that $v_\lambda\geq0$ in $D_\lambda$
for any $\lambda>\lambda_*$ sufficiently close to
$\lambda_*$.
This contradicts the definition of $\lambda_*$.
Thus $\lambda_*=a$.
Since $a\in\R$ is arbitrary, we conclude that
$\partial_yu\geq0$ in $\Omega$.
For any fixed $(x'_0,y_0)\in\Omega$, we therefore have
\[
 u(x'_0,y)\geq u(x'_0,y_0)>0
 \ \ \text{for any}\ \ y\geq y_0.
\]
We further deduce from \eqref{eq-epigraph-domain} that
$\dist((x'_0,y),\Omega^c)\to+\infty$
as $y\to+\infty$, which contradicts the first
alternative in \eqref{eq-deep-alts}.
Consequently, the second alternative holds.
For each $\theta\in(0,1)$, choose $R>0$ such that
$\inf_{\Omega_R}u\geq\theta$.
Then \eqref{eq-domain-bounds} holds with $R_\theta=R+1$.
This completes the proof of Lemma~\ref{lem-domain-bounds}.
\end{proof}

Let $K\subset\R^{d-1}$ be compact and let $R>0$. If
$y>R+\max_{\dist(\cdot,K)\leq R}\varphi$, then \eqref{eq-epigraph-domain}
implies that
\begin{equation}\label{eq-vertical-distance}
 B_R((x',y))\subset\Omega\ \ \text{for any}\ \ x'\in K.
\end{equation}
Choose $R$ with $2d/R^2\leq\min_{[2/(d+2)^2,\,1/2]}f$. 
It follows from \eqref{eq-vertical-distance}
that there exists $x_0\in\Omega$ such that
$B_R(x_0)\Subset\Omega$. Define
\[
 w_0(x):=\tfrac12(1-|x-x_0|^2/R^2)_+^2
 \,\ \ \text{in}\ \ \R^d.
\]
For any $x\in B_R(x_0)$,  we have $-\Delta
w_0\leq2d/R^2$, with $-\Delta w_0>0$ only where $w_0\geq2/(d+2)^2$. 
Since $w_0=0$ and $\nabla w_0=0$ on $\partial B_R(x_0)$,, we obtain that
\begin{equation}\label{eq-general-ball-subsolution}
 -\Delta w_0\leq f(w_0)\,\ \ \text{weakly in}\ \ \Omega.
\end{equation}
Set $M:=1+\max_{[0,1]}|f'|$. For any measurable $h$ with $0\leq h\leq1$, let
$Th$ be the bounded weak solution of
\begin{equation}\label{eq-epigraph-iteration}
 (-\Delta+M)Th=Mh+f(h)\,\ \ \text{in}\ \ \Omega,\ \ 
 Th\equiv0\,\ \ \text{on}\ \ \partial\Omega .
\end{equation}
Let $q,z\in H^1_{\mathrm{loc}}(\overline\Omega)$ be a weak subsolution and a
weak supersolution of \eqref{mainequation}, respectively, with $0\leq q\leq
z\leq1$ in $\Omega$ and traces $q\equiv0$ and $z\geq0$ on $\partial\Omega$.
Since $s\mapsto Ms+f(s)$ is increasing on $[0,1]$, \cite[Lemma~2.1]{JW2026}
implies that $q\leq Tq$, $Tz\leq z$ and $T$ preserves order in the sense that
\begin{equation}\label{eq-epigraph-order-preserving}
 0\leq h_1\leq h_2\leq1
 \ \ \Longrightarrow\ \ 
 Th_1\leq Th_2
 \,\ \ \text{a.e. in}\ \ \Omega
\end{equation}
for any measurable $h_1,h_2$ on $\Omega$. In particular,
\begin{equation}\label{eq-epigraph-ordered-iteration}
 q\leq T^nq\leq T^{n+1}q\leq T^nz\leq
 z\,\ \ \text{a.e.\ in}\ \ \Omega\ \ \text{for any}\ \ n\geq0.
\end{equation}
As in the proof of \cite[Lemma~2.3]{JW2026}, the limit of $T^nq$ solves
\eqref{mainequation}. In particular, monotone iteration from $w_0$ and $1$
yields positive bounded solutions $\underline u$ and $\overline u$ of
\eqref{mainequation} with
\begin{equation}\label{eq-extremal-limits}
 T^nw_0\to \underline u,\ \  T^n1\to \overline
 u\,\ \ \text{in}\ \ C^2_{\mathrm{loc}}(\Omega)\ \ \text{as}\ \ n\to\infty.
\end{equation}

\begin{prop}\label{prop-extremal}
Under the assumptions of Theorem~{\ref{thm-main}}(i), suppose that $f'(0)>0$ or
that $d=2$ and \eqref{eq-planar-quadratic} holds. Then there exist positive
bounded solutions $\underline u,\overline u$ of \eqref{mainequation} such that
each positive bounded solution $u$ satisfies
\begin{equation}\label{eq-extremal-order}
 0<\underline u\leq u\leq \overline u<1\,\ \ \text{in}\ \ \Omega.
\end{equation}
Moreover, for each $w\in\{\underline u,\overline u\}$ we have
\begin{equation}\label{eq-vertical-extremal}
 \partial_yw>0\,\ \ \text{in}\ \ \Omega,\ \ 
 w(x',y)\to1\ \ \text{as}\ \ y\to+\infty,
\end{equation}
where the convergence is locally uniform in $x'$.
\end{prop}

\begin{proof}
Let $u$ be any positive bounded solution of
\eqref{mainequation}.
We first prove \eqref{eq-extremal-order}.
It follows from \eqref{eq-domain-bounds} and
\eqref{eq-vertical-distance} that
$u(x',y)\to1$ locally uniformly in $x'$ as $y\to+\infty$.
Since $w_0\leq1/2$ and $w_0$ has compact support,
we have $w_0(\cdot-te_d)\leq u$ in $\Omega$
for any sufficiently large $t>0$.
The translated function $w_0(\cdot-te_d)$ remains
a subsolution of \eqref{mainequation} as $t$ decreases.
Since $u>0$ on $\partial B_R(x_0+te_d)$,
the strong maximum principle excludes a first contact
for $t>0$.
The comparison therefore holds for any $t\geq0$,
which implies that
\begin{equation}\label{eq-extremal-subsuper}
 w_0(x)\leq u(x+te_d)
 \ \ \text{for any}\ \ x\in\Omega,\ \ t\geq0.
\end{equation}
Taking $t=0$ in \eqref{eq-extremal-subsuper},
we then deduce from \eqref{eq-epigraph-ordered-iteration}
that
\[
 T^nw_0\leq u\leq T^n1
 \,\ \ \text{in}\ \ \Omega
 \ \ \text{for any}\ \ n\geq1.
\]
Letting $n\to\infty$, we then obtain from
\eqref{eq-extremal-limits} that
\eqref{eq-extremal-order} holds.

We next prove the strict vertical monotonicity of
$\underline u$ and $\overline u$.
For any $t>0$, $\underline u(\cdot+te_d)$ is
a supersolution of \eqref{mainequation} in $\Omega$
and satisfies $w_0\leq\underline u(\cdot+te_d)$
in $\Omega$ according to \eqref{eq-extremal-subsuper}.
The zero extension of $\overline u(\cdot-te_d)$
from $\Omega+te_d$ to $\Omega$ is a weak subsolution
of \eqref{mainequation} and takes values in $[0,1]$.
We derive from \eqref{eq-epigraph-ordered-iteration}
and \eqref{eq-extremal-limits} that
\[
 w(x+te_d)\geq w(x)
 \,\ \ \text{in}\ \ \Omega
 \ \ \text{for any}\ \ t>0
\]
for each $w\in\{\underline u,\overline u\}$.
It follows from \eqref{mainequation} that
\[
 \partial_yw\geq0,\ \ 
 -\Delta\partial_yw=f'(w)\partial_yw
 \,\ \ \text{in}\ \ \Omega.
\]
Since $w>0$ in $\Omega$ and $w=0$ on $\partial\Omega$,
the strong maximum principle implies that
$\partial_yw>0$ in $\Omega$.
The convergence in \eqref{eq-vertical-extremal}
has already been established for any positive bounded
solution of \eqref{mainequation}.
This proves \eqref{eq-vertical-extremal}, which completes
the proof of Proposition~\ref{prop-extremal}.
\end{proof}

\begin{proof}[Proof of Theorem~\ref{thm-main}(i)]
Existence follows from \eqref{eq-extremal-limits}.
Under
\eqref{eq-planar-quadratic},
let $d=2$,
$u:=\underline u$ and $v:=\overline u$
be as in Proposition~\ref{prop-extremal}.
Suppose on the contrary that $u\not\equiv v$ in $\Omega$.
Applying the strong maximum principle,
we then deduce from \eqref{eq-extremal-order}
that $u<v$ in $\Omega$.
Define
\[
 G(x')
 :=\tfrac12(|\nabla v|^2-|\nabla u|^2)
 (x',\varphi(x'))
 \ \ \text{for any}\ \ x'\in\R.
\]
We first prove that $G\in L^1(\R)$.
Using the interior and boundary elliptic estimates
\cite[Corollary~8.36]{GT2001},
we then obtain a constant $M>0$ such that
\begin{equation}\label{eq-epigraph-gradient-bound}
 |\nabla u|,|\nabla v|\leq M
 \,\ \ \text{on}\ \ \overline\Omega.
\end{equation}
Since $f(1)=0$, we derive from
\eqref{mainequation} and \eqref{eq-vertical-extremal} that
\begin{equation}\label{eq-epigraph-top-limit}
 (\nabla u,\nabla v)(x',y)\to(0,0)
 \ \ \text{locally uniformly in}\ \ x'
 \ \ \text{as}\ \ y\to+\infty.
\end{equation}
It follows from \eqref{eq-vertical-extremal}
that for any $x'\in\R$ and $0<s<1$,
there exists a unique $h_w(x',s)$ such that
$w(x',h_w(x',s))=s$, where $w=u,v$.
Changing variables, we then obtain that
\begin{equation}\label{eq-epigraph-flux}
 \begin{aligned}
 J(x')
 &:=\int_{\varphi(x')}^{+\infty}
 (v_{x'}v_y-u_{x'}u_y)(x',y)\,dy\\
 &=\int_0^1
 [v_{x'}(x',h_v(x',s))-u_{x'}(x',h_u(x',s))]\,ds.
 \end{aligned}
\end{equation}
We thus obtain $J\in C(\R)$ and $|J|<2M$.
Define $F(s):=\int_0^s f(t)\,dt$ and
$S_w:=(w_{x'}w_y,(w_y^2-w_{x'}^2)/2+F(w))$
for each $w=u,v$.
We calculate from \eqref{mainequation} that
$\diver S_w=w_y(\Delta w+f(w))=0$ in $\Omega$
and $S_w\cdot\nu=\tfrac12|\nabla w|^2\nu_y$
on $\partial\Omega$, where $\nu$ is the outer unit normal.
Applying the Hopf lemma to $v-u$, we then obtain
$\partial_\nu v<\partial_\nu u\leq0$ on $\partial\Omega$,
which implies that
\begin{equation}\label{eq-epigraph-boundary-gap}
 G(x')>0
 \ \ \text{for any}\ \ x'\in\R.
\end{equation}
Applying the divergence theorem, we then obtain that
for any $\eta\in C_c^\infty(\R)$
and $R>\max_{\supp\eta}\varphi$,
\begin{equation*}
 \begin{aligned}
 &\int_{\Omega\cap\{y<R\}}
 (v_{x'}v_y-u_{x'}u_y)\eta'(x')\,dx'\,dy\\
 &=\int_{\Omega\cap\{y<R\}}
 \diver\bigl(\eta(S_v-S_u)\bigr)\,dx'\,dy\\
 &=\int_{\partial\Omega\cap\{y<R\}}
 G(x')\eta(x')\nu_y\,dS
 +\int_{\supp\eta}
 (S_v-S_u)_2(x',R)\eta(x')\,dx'\\
 &=-\int_\R G(x')\eta(x')\,dx'
 +\int_{\supp\eta}
 (S_v-S_u)_2(x',R)\eta(x')\,dx',
 \end{aligned}
\end{equation*}
where the last equality follows from the area formula
\cite[Theorems~3.8 and~3.9]{EG2015}
and $(\cdot)_2$ denotes the second component.
It follows from \eqref{eq-vertical-extremal}
and \eqref{eq-epigraph-top-limit} that
the integral over $\{y=R\}$ tends to zero
as $R\to+\infty$.
Applying dominated convergence, we then deduce from
\eqref{eq-epigraph-gradient-bound}
and \eqref{eq-epigraph-flux} that
\begin{equation*}
 \int_\R J\eta'\,dx'
 =\lim_{R\to+\infty}
 \int_{\Omega\cap\{y<R\}}
 (v_{x'}v_y-u_{x'}u_y)\eta'(x')\,dx'\,dy
 =-\int_\R G\eta\,dx'.
\end{equation*}
Thus $J'=G$ in $\R$.
Since $G>0$ according to
\eqref{eq-epigraph-boundary-gap} and $|J|\leq2M$,
we obtain that $J$ is strictly increasing and
\begin{equation}\label{eq-epigraph-integrable-gap}
 \int_\R G(x')\,dx'\leq4M.
\end{equation}

We now establish two distinct positive solutions
of \eqref{mainequation}.
Reflecting the horizontal coordinate if necessary,
we then obtain a constant $j>0$ such that
$|J(x')|\geq j$ for any sufficiently large $x'>0$.
Choose $x_n'$ as a maximum point of $\varphi$ on $[n,2n]$,
so that $x_n'\to+\infty$ as $n\to\infty$.
Choose $\delta\in(0,1/2)$ such that $4M\delta<j/2$.
We deduce from \eqref{eq-epigraph-gradient-bound}
and \eqref{eq-epigraph-flux} that
for any sufficiently large $n>0$,
\begin{equation*}
 \int_\delta^{1-\delta}
 |v_{x'}(x_n',h_v(x_n',s))
 -u_{x'}(x_n',h_u(x_n',s))|\,ds
 \geq |J(x_n')|-4M\delta>j/2.
\end{equation*}
Hence there exists $s_n\in[\delta,1-\delta]$ such that
\begin{equation}\label{eq-epigraph-derivative-gap}
 |v_{x'}(x_n',h_v(x_n',s_n))
 -u_{x'}(x_n',h_u(x_n',s_n))|\geq j/2
 \ \ \text{for any sufficiently large}\ \ n>0.
\end{equation}
Define
\begin{equation*}
 \bigl(u_n(x',y),v_n(x',y)\bigr)
 :=\Bigl(
 u\bigl(x'+x_n',y+h_u(x_n',s_n)\bigr),
 v\bigl(x'+x_n',y+h_u(x_n',s_n)\bigr)
 \Bigr)
 \,\ \ \text{in}\ \ \Omega_n,
\end{equation*}
where $\Omega_n:=\Omega-(x_n',h_u(x_n',s_n))$.
It follows from \eqref{eq-epigraph-gradient-bound} that
\begin{equation*}
 \dist\bigl((x_n',h_w(x_n',s_n)),\partial\Omega\bigr)
 \geq s_n/M\geq4r
 \ \ \text{for each}\ \ w=u,v,
 \ \ r:=\delta/(4M).
\end{equation*}
Define
\[
 W_n(x',y)
 :=v\bigl(x'+x_n',y+h_v(x_n',s_n)\bigr)-u_n(x',y)
 \,\ \ \text{in}\ \ \overline B_{3r}.
\]
It follows from Proposition~\ref{prop-extremal} that
$h_v(x_n',s_n)\leq h_u(x_n',s_n)$.
We derive from \eqref{eq-vertical-extremal}
and \eqref{eq-epigraph-derivative-gap} that
for any sufficiently large $n>0$,
\begin{equation}\label{eq-epigraph-translated-difference}
 W_n\leq v_n-u_n
 \,\ \ \text{in}\ \ B_{3r},
 \ \ W_n(0)=0,
 \ \ |\partial_{x'}W_n(0)|\geq j/2.
\end{equation}
It follows from \cite[Proposition~A.4]{BG2025} that
for any $0<\alpha<\gamma$,
$\Omega_n\to\Omega_\infty$ locally in $C^{2,\alpha}$
as $n\to\infty$ along a subsequence.
Let $D$ be the connected component of $\Omega_\infty$
containing $\overline B_{2r}$.
Passing to a further subsequence, we then obtain that
\begin{equation*}
 (u_n,v_n)\to(u_\infty,v_\infty)
 \,\ \ \text{in}\ \ C^2_{\mathrm{loc}}(D)^2
 \ \ \text{and in}\ \ C^1_{\mathrm{loc}}(\overline D)^2
 \ \ \text{as}\ \ n\to\infty.
\end{equation*}
Moreover, $W_n\to W_\infty$ in $C^1(\overline B_r)$
as $n\to\infty$ along a further subsequence.
The limits $u_\infty,v_\infty$ are solutions
of \eqref{mainequation} in $D$,
with $u_\infty(0)\geq\delta$.
If $u_\infty\equiv v_\infty$ in $D$,
we deduce from \eqref{eq-epigraph-translated-difference}
that $W_\infty\leq W_\infty(0)=0$ in $B_r$,
and hence
$0=|\partial_{x'}W_\infty(0)|\geq j/2$,
a contradiction.
Applying the strong maximum principle,
we then obtain that
$0<u_\infty<v_\infty\leq1$ in $D$.

We finally derive a contradiction from the geometry of $D$.
The outer normals to the translated epigraphs
have nonpositive vertical component.
Local $C^1$ convergence of the boundaries therefore implies that
$\nu_y\leq0$ on $\partial D$.
If $\nu_y<0$ at some point of $\partial D$,
we can represent the boundary near that point by
$y=\varphi_\infty(x')$ on a bounded interval $I_0$
with closure in the coordinate neighborhood.
Using the $C^1$ convergence up to this boundary,
we then deduce from \eqref{eq-epigraph-integrable-gap} that
\begin{equation*}
 0=\lim_{n\to\infty}\int_{I_0}G(x_n'+x')\,dx'
 =\tfrac12\int_{I_0}
 (|\nabla v_\infty|^2-|\nabla u_\infty|^2)
 (x',\varphi_\infty(x'))\,dx'>0,
\end{equation*}
where the strict inequality follows from the Hopf lemma.
Hence $\nu_y\equiv0$ on $\partial D$,
so the boundary consists of vertical lines.
Applying the uniqueness results in the plane
\cite[Section~8.4]{BG2026}
and in the half-plane \cite[Theorem~1.5]{BCNHalf1997},
and using $(w')^2=2(F(1)-F(w))$
for the half-plane profile,
we then obtain that
\begin{equation*}
 D=(a,b)\times\R,
 \ \ -\infty<a<0<b<+\infty.
\end{equation*}
For any sufficiently small $\rho>0$
and any sufficiently large $n$,
the points $(a-\rho,1)$ and $(b+\rho,1)$
lie outside $\Omega_n$.
Hence
\begin{equation*}
 \min\{\varphi(x_n'+a-\rho),\varphi(x_n'+b+\rho)\}
 \geq h_u(x_n',s_n)+1>\varphi(x_n').
\end{equation*}
At least one of $x_n'+a-\rho$ and $x_n'+b+\rho$
belongs to $[n,2n]$ for any sufficiently large $n$,
since $\max\{x_n'-n,2n-x_n'\}\geq n/2$.
This contradicts the choice of $x_n'$.
We thus obtain that $u\equiv v$ in $\Omega$.
It follows from \eqref{eq-extremal-order} that
the positive bounded solution of \eqref{mainequation}
is unique.
This completes the proof of Theorem~\ref{thm-main}(i).
\end{proof}

\section{Nonuniqueness for \texorpdfstring{$N \ge 3$}{N>=3}}\label{sec-high}

{
The purpose of this section is to prove
Theorem~\ref{thm-main}(ii) and Theorem~\ref{thm-extension}.
}

\begin{lem}\label{lem-transverse-states}
We have the following assertions.
\begin{enumerate}
 \item[(i)] Let $f$ satisfy the assumptions of Theorem~{\ref{thm-main}}(ii).
 There exist an interval $I:=(-L,L)$ and a positive function $B\in H^1_0(I)$
 such that $0$ and $B$ are the only nonnegative global minimizers of $\mathcal
 J(v):=\int_I(\tfrac12|v'|^2-F(|v|))dx_1$ on $H^1_0(I)$, where $F(s):=\int_0^s
 f(t)dt$. Moreover, the first Dirichlet eigenvalues of
 $-\partial_{x_1}^2-f'(0)$ and $-\partial_{x_1}^2-f'(B)$ on $I$ are positive.
 \item[(ii)] For $I$ and $B$ in \textup{(i)}, $-W_{x_1x_1}-W_{ss}=f(W)$ in
 $I\times\R$ admits a solution $W$ with $W\equiv0$ on $\partial I\times\R$,
 increasing strictly in $s$ from $0$ to $B$. Moreover, there exist $C,\mu>0$
 such that
\begin{equation}\label{eq-strip-tails}
\begin{aligned}
& \norm{W}_{C^2(\overline I\times[s-1,s+1])}\leq Ce^{\mu s}\ \ \text{for
any}\ \ s\leq0,\\
 & \norm{W-B}_{C^2(\overline I\times[s-1,s+1])}\leq Ce^{-\mu s}\ \ \text{for
 any}\ \ s\geq0.
 \end{aligned}
 \end{equation}

{
 \item[(iii)] {
 Suppose $s_0>0$ and $g\in C^{1,\gamma}([0,s_0])$ satisfy
 $g(0)=0$ and $g(s)>0$ for any sufficiently small $s>0$.
 Then there exists a reaction $f$ satisfying \eqref{eq-reaction}
 with $f'(0)=g'(0)\geq0$ in place of $f'(0)>0$ and
 $f=g$ on $[0,\delta]$ for some $\delta\in(0,s_0)$ such that
 the conclusions of \textup{(i)} and \textup{(ii)} hold for this $f$.
 }
}

\end{enumerate}
\end{lem}

\begin{proof}
\textup{(i)} For any $0<m<1$, let $b_m$ be the even solution of 
\begin{equation*}
-b_m''=f(b_m)\,\ \ \text{in}\ \ (-\ell(m),\ell(m)),\ \ 
b_m(0)=m,\ \ b_m'(0)=0,
\end{equation*}
where $\ell(m)>0$ is the smallest positive zero point of $b_m$.
We first prove that $\ell$ is strictly convex.
Since
$(b_m')^2=2(F(m)-F(b_m))$, we have $b_m'(\ell(m))=-\sqrt{2F(m)}<0$. One can
further check that $\ell\in C^1(0,1)$. Set for any $0\leq t\leq1$,
\begin{equation}\label{eq-strip-time-integrand}
 G(m,t):=2\int_t^1\sigma q(m\sigma)d\sigma,
\end{equation}
so that
\begin{equation}\label{eq-strip-time-map}
 \ell(m)=\int_0^m\frac{dr}{\sqrt{2(F(m)-F(r))}}=\int_0^1 G(m,t)^{-1/2}dt.
\end{equation}
Set $q(0):=f'(0)$. Since $q(s)>q(0)$ for some $s\in(0,1)$  and $q$ is concave,
$q$ is strictly increasing near zero. It follows from \eqref{eq-strip-time-integrand} that for any $0<m_1<m_2<1$,
\begin{equation*}
 \partial_t[G(m_2,t)-G(m_1,t)]=-2t[q(m_2t)-q(m_1t)]<0 \ \ \text{for any
 sufficiently small}\ \ t>0.
\end{equation*}
Using the concavity of $G(\cdot,t)$ and the strict convexity
and monotonic decrease of $r\mapsto r^{-1/2}$,
we then deduce from \eqref{eq-strip-time-map}
that $\ell$ is strictly convex.

We next choose $I$ and $B$ and identify
the nonnegative minimizers of $\mathcal J$.
It follows from \eqref{eq-strip-time-map} that  $\ell(m)\to\ell(0)$ as
$m\to0^+$, where $\ell(0):=\pi/(2\sqrt{f'(0)})$ if $f'(0)>0$ and
$\ell(0):=+\infty$ if $f'(0)=0$. Since  $q$ is strictly increasing near zero,
we have $\ell(m)<\ell(0)$ for any sufficiently small $m>0$ when $f'(0)>0$. 
Since $F(1)-F(r)\sim-f'(1)(1-r)^2/2$ as $r\to1^-$,
we derive from \eqref{eq-strip-time-map} that
\begin{equation*}
 \ell(m)\geq\int_0^m[2(F(1)-F(r))]^{-1/2}dr\to+\infty\ \ \text{as}\ \ m\to1^-.
\end{equation*}
Thus $\ell$ has a unique minimum point $m_0\in(0,1)$. Noting that $(b_m')^2=2(F(m)-F(b_m))$, we obtain that
\begin{equation}\label{eq-strip-profile-energy}
 \begin{aligned}
 \mathcal
 E(m):=&\int_{-\ell(m)}^{\ell(m)}\bigl(\tfrac12|b_m'|^2-F(b_m)\bigr)dx_1\\
 =&2\int_0^m\sqrt{2(F(m)-F(r))}\,dr-2F(m)\ell(m)\\
 =&2\int_0^m f(r)(\ell(r)-\ell(m))dr.
 \end{aligned}
\end{equation}
Since $\ell$ is strictly convex, we deduce from
\eqref{eq-strip-profile-energy} that
$\mathcal E(m)>0$ for any $0<m\leq m_0$ and
\begin{equation*}
 \begin{gathered}
 \mathcal E'(m)=-2F(m)\ell'(m)<0\ \ (m_0<m<1),\\
 \mathcal
 E(m)\leq2\sqrt{2F(1)}-2F(m)\ell(m)\to-\infty\ \ \text{as}\ \ m\to1^-.
 \end{gathered}
\end{equation*}
Thus $\mathcal E$ has a unique zero $m_*>m_0$. If $\ell(0)\leq\ell(m_*)$,
the strict convexity of $l$ implies that $\ell(r)<\ell(m_*)$ for any $0<r<m_*$, contrary to
$\mathcal E(m_*)=0$ in \eqref{eq-strip-profile-energy}. We therefore obtain
that
\begin{equation}\label{eq-strip-selected-interval}
 \begin{gathered}
 0<L:=\ell(m_*)<\ell(0),\ \  I:=(-L,L),\\
 \ell'(m)>0\ \ \text{for any}\ \ m_*\leq m<1,\ \ 
 \{m\in(0,1):\ell(m)=L,\ \ \mathcal E(m)\leq0\}=\{m_*\}.
 \end{gathered}
\end{equation}
Applying the direct method in the calculus of variations
\cite[Theorem~2.3]{Rindler2018}
and using the fact that $\mathcal J(|v|)=\mathcal J(v)$,
we then obtain a nonnegative minimizer of $\mathcal J$
on $H^1_0(I)$. Each nonzero nonnegative minimizer is $b_m|_I$ with
$0<m<1$, $\ell(m)=L$ and $\mathcal E(m)\leq0$. It follows from \eqref{eq-strip-selected-interval}
that $m=m_*$. Since $\mathcal
J(B)=\mathcal E(m_*)=0$ for $B:=b_{m_*}|_I$, we obtain that
\begin{equation}\label{eq-strip-energy-minimum}
\inf_{H^1_0(I)}\mathcal J=0,\ \ \{v\in H^1_0(I):v\geq0,\ \ \mathcal
J(v)=0\}=\{0,B\}.
\end{equation}

We next prove that both first Dirichlet eigenvalues
in part \textup{(i)} are positive.
Differentiating $-b_m^{\prime\prime}=f(b_m)$ in $(-\ell(m),\ell(m))$ and
$b_m(\pm\ell(m))=0$ with respect to $m$ at $m=m_*$,
we then obtain from
\eqref{eq-strip-selected-interval} that
\begin{equation*}
 A_Bz=0\,\ \ \text{in}\ \ I,\ \  z(\pm L)=-B'(L)\ell'(m_*)>0,
\end{equation*}
where $A_B:=-\partial_{x_1}^2-f'(B)$ and $z:=\partial_m b_m|_{m=m_*}$.
Minimality of $B$ implies that $\lambda_1(A_B;I)\geq0$. 
Let $\phi$ be a positive first Dirichlet eigenfunction
of $A_B$ on $I$.
The Hopf lemma implies that $\partial_\nu\phi<0$
on $\partial I$.
Since $z(\pm L)>0$, multiplying $A_Bz=0$ by $\phi$
and integrating by parts over $I$, we obtain that
\begin{equation*}
 \lambda_1(A_B;I)\int_I z\phi\,dx_1=-\sum_{x_1\in\partial
 I}z(x_1)\partial_\nu\phi(x_1)>0.
\end{equation*}
Since it follows from \eqref{eq-strip-selected-interval} that
$L<\ell(0)=\pi/(2\sqrt{f'(0)})$ when $f'(0)>0$, we also have
\begin{equation}\label{eq-strip-stability}
 \lambda_1(-\partial_{x_1}^2-f'(0);I)=\frac{\pi^2}{4L^2}-f'(0)>0,\ \ 
 \lambda_1(A_B;I)>0.
\end{equation}
This proves part \textup{(i)}.

\textup{(ii)} We first construct the solution $W$ in part \textup{(ii)}.
Since $\mathcal J(v)=\mathcal J(|v|)$,
we deduce from \eqref{eq-strip-energy-minimum}
that $0$ and $\pm B$ are the only global minimizers
of $\mathcal J$.
Thus $0$ is isolated and $B$ is the unique positive minimizer.
Applying \cite[Theorem~1.1]{DRR2025} to the odd extension of $f$,
we then obtain $W\in C^{2,\gamma}_{\mathrm{loc}}(\overline I\times\R)$
such that
\begin{equation}\label{eq-strip-front}
 \begin{gathered}
 -W_{x_1x_1}-W_{ss}=f(W),\ \  0<W<B,\ \ 
 W_s>0\,\ \ \text{in}\ \ I\times\R,\ \ 
 W\equiv0\,\ \ \text{on}\ \ \partial I\times\R,\\
 W(\cdot,s)\to0\ \ \text{as}\ \ s\to-\infty,\ \  W(\cdot,s)\to
 B\ \ \text{as}\ \ s\to+\infty \,\ \ \text{uniformly on}\ \ \overline I.
 \end{gathered}
\end{equation}

We next prove \eqref{eq-strip-tails}.
Define for each $A\in\{0,B\}$,
\[
 Z(x_1,t):=
 \begin{cases}
  W(x_1,-t), & A=0,\\
  B(x_1)-W(x_1,t), & A=B,
 \end{cases}
 \,\ \ \text{in}\ \ I\times\R,
\]
and
\[
 c(x_1,t):=
 \begin{cases}
  \displaystyle\int_0^1 f'(\theta Z(x_1,t))\,d\theta,
  & A=0,\\
  \displaystyle\int_0^1
  f'(B(x_1)-\theta Z(x_1,t))\,d\theta,
  & A=B,
 \end{cases}
 \,\ \ \text{in}\ \ I\times\R.
\]
Since $-A''=f(A)$ in $I$, we deduce from
\eqref{eq-strip-front} that
\begin{equation*}
 \bigl(-\partial_{x_1x_1}-\partial_{tt}-c\bigr)Z=0
 \,\ \ \text{in}\ \ I\times\R,
 \ \ Z=0
 \,\ \ \text{on}\ \ \partial I\times\R.
\end{equation*}
The limits in \eqref{eq-strip-front} imply that $Z(\cdot,t)\to0$ and
$c(\cdot,t)\to f'(A)$ uniformly on $\overline I$ as $t\to+\infty$. Moreover,
we obtain from \eqref{eq-strip-stability} that there exist  positive functions
$\psi_0$ and $\psi_B$ such that
\begin{equation}\label{eq-strip-positive-tests}
 (-\partial_{x_1}^2-f'(0))\psi_0=A_B\psi_B=1\,\ \ \text{in}\ \ I,\ \ 
 \psi_0\equiv\psi_B\equiv1\,\ \ \text{on}\ \ \partial I.
\end{equation}
Choose $\mu>0$ such that $\mu^2\max\{\|\psi_0\|_\infty,\|\psi_B\|_\infty\}\leq1/4$.
Since $c(\cdot,t)\to f'(A)$ uniformly on $\overline I$
as $t\to+\infty$, we deduce from
\eqref{eq-strip-positive-tests} that for each $A\in\{0,B\}$
and any sufficiently large $T$,
\begin{equation*}
 (-\partial_{x_1x_1}-\partial_{tt}-c)(e^{-\mu t}\psi_A)
 =e^{-\mu t}\bigl[1+(f'(A)-c-\mu^2)\psi_A\bigr]\geq\tfrac12e^{-\mu t}
 \,\ \ \text{in}\ \ I\times(T,+\infty).
\end{equation*}
Since $(-\partial_{x_1x_1}-\partial_{tt}-c)\psi_A>0$
in $I\times(T,+\infty)$, we compare $Z$ with
$Ce^{-\mu t}\psi_A$ for any sufficiently large $C$
and obtain that
\begin{equation}\label{eq-strip-tail-comparison}
 0\leq Z(x_1,t)\leq Ce^{-\mu t}\psi_A(x_1)
 \ \ \text{for any}\ \ x_1\in\overline I,\ \ t\geq T.
\end{equation}
Applying interior and boundary Schauder estimates
\cite[Corollaries~6.3 and~6.7]{GT2001},
we then deduce from \eqref{eq-strip-tail-comparison}
that \eqref{eq-strip-tails} holds.
{
This proves part \textup{(ii)}.
}

{
\textup{(iii)} Take $f(s):=(1+g'(0))(s+s^2-2s^3)$. Let $\ell$, $\mathcal E$, $m_0$ and $m_*$ be as in
part \textup{(i)} for this reaction.
Choose $m_{\mathrm L}$ and $m_{\mathrm R}$ such that
$0<m_0<m_{\mathrm L}<m_*<m_{\mathrm R}<1$. Since $f'(0)=1+g'(0)$ and $\ell(0)>\ell(m_*)$, there
exists $s_1\in(0,\min\{m_{\mathrm L},s_0\})$ such that $0<g(s)\leq f(s)$ for
any $0<s\leq s_1$ and $\inf_{0<m\leq s_1}\ell(m)>\ell(m_*)$. Let $\chi\in
C^\infty([0,+\infty);[0,1])$ satisfy that $\chi\equiv 1$ on $[0,1]$ and
$\chi\equiv 0$ on $[2,+\infty)$.  Set for any $0<\delta<s_1/4$, 
\begin{equation*}
f_\delta(s):=f(s)+\chi(s/\delta)(g(s)-f(s))\,\ \ \text{on}\ \ [0,s_1],\ \ 
f_\delta(s):=f(s)\,\ \ \text{on}\ \ (s_1,+\infty).
\end{equation*} 
Then $f_\delta$ satisfies \eqref{eq-reaction} with $f_\delta'(0)=g'(0)\geq0$
in place of $f_\delta'(0)>0$ and $f_\delta=g$ on $[0,\delta]$.

Set $F_\delta(s):=\int_0^s f_\delta(t)dt$. For any $s_1\leq m<1$, we have
$F_\delta(m)-F_\delta(r)=F(m)-F(r)$ for any $2\delta\leq r\leq m$. For any
$0\leq r\leq2\delta$, $F_\delta(m)-F_\delta(r)$ and $F(m)-F(r)$ are bounded
below by $\int_{s_1/2}^{s_1}f(t)dt>0$ and differ by at most
$\int_0^{2\delta}|f_\delta(t)-f(t)|dt\leq C\delta^2$. For $f_\delta$, define
$\ell_\delta$ by \eqref{eq-strip-time-map} and $\mathcal E_\delta$ by
\eqref{eq-strip-profile-energy}. We then have
\begin{equation}\label{eq-strip-extension-convergence}
 \sup_{s_1\leq
 m<1}\bigl(|\ell_\delta(m)-\ell(m)|+|\ell_\delta'(m)-\ell'(m)|\bigr)\leq
 C\delta^3,\ \ 
 \sup_{s_1\leq m\leq m_{\mathrm R}}|\mathcal E_\delta(m)-\mathcal E(m)|\leq
 C\delta^2 .
\end{equation}
Since $\inf_{[m_{\mathrm L},1)}\ell'=\ell'(m_{\mathrm L})>0$,
\eqref{eq-strip-extension-convergence} gives that for any sufficiently small
$\delta>0$,
\begin{equation*}
 \ell_\delta'>0\,\ \ \text{on}\ \ [m_{\mathrm L},1),\ \ 
 \min_{[s_1,m_{\mathrm L}]}\mathcal E_\delta>0>\mathcal E_\delta(m_{\mathrm
 R}).
\end{equation*}
Since $\mathcal E_\delta'=-2F_\delta\ell_\delta'$, the unique zero $m_\delta$
in $[s_1,1)$ belongs to $(m_{\mathrm L},m_{\mathrm R})$. Uniform convergence
in \eqref{eq-strip-extension-convergence} gives $m_\delta\to m_*$ and
$\ell_\delta(m_\delta)\to\ell(m_*)$ as $\delta\to0^+$. For any sufficiently
small $\delta>0$, $f_\delta\leq f$ and \eqref{eq-strip-time-map} then give
that
\begin{equation*}
 \ell_\delta(m_\delta)<\inf_{0<m\leq s_1}\ell(m)\leq\inf_{0<m\leq
 s_1}\ell_\delta(m)\leq\ell_\delta(0).
\end{equation*}
Thus \eqref{eq-strip-selected-interval} holds with $m_\delta$ and
$L=\ell_\delta(m_\delta)$. The rest of \textup{(i)} and \textup{(ii)} now
apply to $f_\delta$. Fixing $\delta>0$ sufficiently small,
we then take $f=f_\delta$ in part \textup{(iii)}.
This completes the proof of
Lemma~\ref{lem-transverse-states}.
}
\end{proof}

\begin{lem}\label{lem-strip-cap}
{
Let $f$, $I$ and $B$ be as in
Lemma~\ref{lem-transverse-states}(i) or (iii).
}
Then there exist positive constants $\kappa,k,C$ and $P\in
C^{2,\gamma}_{\mathrm{loc}}(\overline I\times\R)$ satisfying $B\leq P\leq2$
such that
\begin{equation}\label{eq-epi-cap-bounds}
 \begin{gathered}
 P_t\geq0,\ \ 
 -P_{x_1x_1}-P_{tt}-f(P)\geq\kappa|P_{tt}|\,\ \ \text{in}\ \ I\times\R,\ \ 
 P>1\,\ \ \text{on}\ \ \overline I\times[0,+\infty),\\
 0\leq P-B\leq Ce^{kt}\,\ \ \text{on}\ \ \overline I\times(-\infty,0].
 \end{gathered}
\end{equation}
\end{lem}

\begin{proof}
We first construct $P$ with $P_t\geq0$
and determine its limits as $t\to\pm\infty$.
For $\psi_B$ in \eqref{eq-strip-positive-tests},
choose $0<c<-f(2)$ and $k>0$ such that
$k^2\norm{\psi_B}_{L^\infty(I)}<1/4$, and define
\[
\zeta(t):=(1+e^{-kt})^{-1}. 
\]
Since $-B''=f(B)$ in $I$ and $f(2)+c<0$, $B$ and
the constant $2$ are a subsolution and a supersolution, respectively, of 
\begin{equation}\label{eq-epi-cap-forced}
 -P_{x_1x_1}-P_{tt}=f(P)+c\zeta(t)\,\ \ \text{in}\ \ I\times\R,\ \ 
 P=2\zeta(t)\,\ \ \text{on}\ \ \partial I\times\R.
\end{equation}
Applying the descending iteration from $2$
as in \eqref{eq-epigraph-iteration}--\eqref{eq-extremal-limits},
we then obtain a solution $B\leq P\leq2$
of \eqref{eq-epi-cap-forced} such that
for any subsolution $v$ of \eqref{eq-epi-cap-forced}
satisfying $B\leq v\leq2$ in $I\times\R$,
\[
 v\leq P
 \,\ \ \text{in}\ \ I\times\R. 
\]

Since $\zeta(t-h)\leq\zeta(t)$ for any $h>0$, $P(x_1,t-h)$ is such a
subsolution. Hence $P_t\geq0$ in $I\times\R$, so for each $x_1\in\overline I$
we may define
\begin{equation*}
 P_{-\infty}(x_1):=\lim_{t\to-\infty}P(x_1,t),\ \ 
 P_{+\infty}(x_1):=\lim_{t\to+\infty}P(x_1,t).
\end{equation*}
By the interior and boundary elliptic estimates
\cite[Corollary~8.36]{GT2001}, we obtain that
$P(\cdot,t)\to P_{\pm\infty}$ uniformly on $\overline I$
as $t\to\pm\infty$. Since $\zeta(t)\to0$ as $t\to-\infty$
and $\zeta(t)\to1$ as $t\to+\infty$,
we deduce from \eqref{eq-epi-cap-forced} that
\begin{align}
 -P_{-\infty}''=f(P_{-\infty})\,\ \ \text{in}\ \ I, \ \ 
 P_{-\infty}\equiv0\,\ \ \text{on}\ \ \partial I,\notag\\
  -P_{+\infty}''=f(P_{+\infty})+c\,\ \ \text{in}\ \ I, \ \ 
  P_{+\infty}\equiv2\,\ \ \text{on}\ \ \partial I\label{eq-equation-Pinf}.
\end{align}
Since $P_{-\infty}\geq B$, its maximum $m$ satisfies $m\geq B(0)=m_*$. The
maximum principle and ODE uniqueness imply that $m<1$. It follows from \eqref{eq-strip-time-map}
and \eqref{eq-strip-selected-interval} that
$\ell(m)=L=\ell(m_*)$ and $\ell'>0$ on $[m_*,1)$,
so $P_{-\infty}\equiv B$. If $P_{+\infty}$ has an interior minimum with value at most $1$,
we obtain from \eqref{eq-equation-Pinf} that
$0\geq-P_{+\infty}''=f(P_{+\infty})+c>0$,
a contradiction.
Thus $\min_{\overline I}P_{+\infty}>1$.

We next prove the differential inequality
and exponential estimate in \eqref{eq-epi-cap-bounds}.
Set $Z(x_1,t):=P(x_1,-t)-B(x_1)$ and $c_P(x_1,t):=\int_0^1 f'(B(x_1)+\theta
Z(x_1,t))\,d\theta$. We deduce from \eqref{eq-epi-cap-forced} that
\begin{equation*}
 (-\partial_{x_1x_1}-\partial_{tt}-c_P)Z=c\zeta(-t)
 \,\ \ \text{in}\ \ I\times\R,\ \ 
 Z=2\zeta(-t)\,\ \ \text{on}\ \ \partial I\times\R.
\end{equation*}
Since $P_{-\infty}\equiv B$, $Z(\cdot,t)\to0$ and $c_P(\cdot,t)\to f'(B)$
uniformly on $\overline I$ as $t\to+\infty$. Comparing $Z$ with $Ce^{-kt}\psi_B$
as in \eqref{eq-strip-tail-comparison}
and applying boundary Schauder estimates,
we then obtain that
\begin{equation*}
 0\leq P(x_1,t)-B(x_1)+|P_{tt}(x_1,t)|\leq Ce^{kt}
 \ \  (x_1\in\overline I,\ \ t\leq0).
\end{equation*}
For any $t\geq0$, boundary estimates for \eqref{eq-epi-cap-forced} imply that
$|P_{tt}|\leq C$. Since $\zeta(t)\geq\tfrac12\min\{e^{kt},1\}$, we obtain from
\eqref{eq-epi-cap-forced} that
\begin{equation*}
 -P_{x_1x_1}-P_{tt}-f(P)=c\zeta(t)\geq\kappa|P_{tt}|
 \,\ \ \text{in}\ \ I\times\R
\end{equation*}
for some $\kappa>0$. Since $P(\cdot,t)\to P_{+\infty}$ uniformly
on $\overline I$ as $t\to+\infty$
and $\min_{\overline I}P_{+\infty}>1$,
we choose $T>0$ such that $P>1$
on $\overline I\times[T,+\infty)$.
Then $P(x_1,t+T)$ satisfies \eqref{eq-epi-cap-bounds}
after increasing $C$ in the exponential estimate,
which completes the proof of Lemma~\ref{lem-strip-cap}.
\end{proof}

\begin{proof}[{Proof of Theorem~{\ref{thm-main}}(ii)}]
We first construct a positive function to perturb
the front $W$ in dimension three.
Write $x=(x',y)$ with $x'=(x_1,x_2)\in\R^2$. Take $B,W$ from
Lemma~\ref{lem-transverse-states}(i)--(ii) and $\chi\in C^\infty(\R;[0,1])$
with $\chi\equiv0$ on $(-\infty,-1]$ and $\chi\equiv1$ on $[1,+\infty)$. Set
$\Psi(x_1,s):=(1-\chi(s))\psi_0(x_1)+\chi(s)\psi_B(x_1)$. Then
$\inf_{\overline I\times\R}\Psi\geq\min_{\overline I}\min\{\psi_0,\psi_B\}>0$,
and \eqref{eq-strip-positive-tests} implies that
\begin{equation*}
 (-\partial_{x_1}^2-\partial_s^2-f'(W))\Psi=
 \begin{cases}
 1,&x_1=\pm L,\\
 1+[f'(0)-f'(W)]\psi_0,&s\leq-1,\\
 1+[f'(B)-f'(W)]\psi_B,&s\geq1.
 \end{cases}
\end{equation*}
Since the displayed expression equals $1$ on the sides and tends uniformly to
$1$ as $|s|\to\infty$ by \eqref{eq-strip-tails}, the set where
$(-\partial_{x_1}^2-\partial_s^2-f'(W))\Psi\leq1/2$ is a compact subset of
$I\times\R$. Since $W_s>0$ there by \eqref{eq-strip-front}, we can take
$C_0>0$ sufficiently large that
\begin{equation}\label{eq-epi-front-perturbation}
 (-\partial_{x_1}^2-\partial_s^2-f'(W))\Psi
 \geq \tfrac14-C_0W_s\,\ \ \text{in}\ \ I\times\R.
\end{equation}

We next construct perturbations of $W$ for the upper
and lower solutions of \eqref{mainequation}.
Fix $K>2(C_0+1)$ and set for any $R\geq4$
\begin{equation}\label{eq-epi-front-parameters}
 h(x_2):=\sqrt{1+x_2^2/R^2},\ \ 
 \omega(x_2):=h''(x_2)/K=\frac{1}{KR^2h(x_2)^3}.
\end{equation}
Set $a(x_2):=(1+h'(x_2)^2)^{-1/2}$. For each $j\in\{-1,1\}$, define
\begin{equation}\label{eq-epi-modulated-fronts}
 s_j(x_2,y):=a(x_2)(y-jh(x_2)).
\end{equation}
We then obtain uniformly for any $R\geq4$ and $x_2\in\R$ that
\begin{equation}\label{eq-epi-modulation-bounds}
 |h'|\leq R^{-1},\ \ \tfrac12\leq a\leq1,\ \ |\omega'|+R|\omega''|\leq
 C\omega/R,\ \ 
 |a'h'|+|a''/a|\leq C\omega/R^2,
\end{equation}
where $C$ may depend on $K$. Moreover, we have
\begin{equation}\label{eq-epi-phase-identities}
 |\nabla_{x_2,y}s_j|^2-1=-2ja'h's_j+(a'/a)^2s_j^2,\ \ 
 \Delta_{x_2,y}s_j+jah''=-2ja'h'+(a''/a)s_j.
\end{equation}
Define
\begin{equation}\label{eq-epi-front-barriers}
 V_j(x',y):=W(x_1,s_j(x_2,y))+j\omega(x_2)\Psi(x_1,s_j(x_2,y)).
\end{equation}
Since $\Psi_s\equiv\Psi_{ss}\equiv0$ for any $|s|>1$,
we obtain from Taylor's formula that in $I\times\R^2$ for $j=1$ and in $\{V_{-1}>0\}$ for $j=-1$,
\begin{equation*}
\begin{aligned}
&\Bigl|j[-\Delta
V_j-f(V_j)]-\omega\Big[aKW_s(x_1,s_j)
-\Psi_{x_1x_1}(x_1,s_j)-\Psi_{ss}(x_1,s_j)\\
&\ \ -f'(W(x_1,s_j))\Psi(x_1,s_j)\Bigr]\Big|\\
\leq&CR^{-2}\omega\big((1+|s_j|)|W_s(x_1,s_j)|
+(1+s_j^2)|W_{ss}(x_1,s_j)|\big)+C(\omega R^{-2}+\omega^2+\omega^{1+\gamma})\\
\leq&C\omega(R^{-2}+\omega^\gamma).
\end{aligned}
\end{equation*}
where we used \eqref{eq-epi-modulation-bounds} and
\eqref{eq-epi-phase-identities} in the first inequality and
\eqref{eq-strip-tails} in the second. It follows from \eqref{eq-epi-front-perturbation}
that for any sufficiently large $R>0$, 
\begin{equation}\label{eq-epi-strict-barriers}
 j[-\Delta V_j-f(V_j)]\geq
 \omega[\tfrac14+(aK-C_0)W_s(x_1,s_j)]
 -C\omega(R^{-2}+\omega^\gamma)\geq\tfrac18\omega
\end{equation}
in $I\times\R^2$ for $j=1$ and in $\{V_{-1}>0\}$ for $j=-1$.

Let $P$ and $\kappa$ be as in
Lemma~\ref{lem-strip-cap}.
Define for any $H_0>1$,
\begin{equation}\label{eq-Hx2}
H(x_2):=H_0+4h(x_2),
\end{equation} 
and $P_H(x',y):=\min\{1,P(x_1,y-H(x_2))\}$.
We next prove that $P_H$ is an upper solution
of \eqref{mainequation} and compare it with $V_1$.
It then follows
from \eqref{eq-epi-front-parameters} and \eqref{eq-epi-modulation-bounds} that
$H''=4h''>0$ and $|H'|^2\leq16/R^2\leq\kappa$ for any sufficiently large
$R>0$. We deduce from \eqref{eq-epi-cap-bounds} that
\begin{equation*}
 \begin{aligned}
 &-\Delta[P(x_1,y-H(x_2))]-f(P(x_1,y-H(x_2)))\\
 &\geq(\kappa-|H'(x_2)|^2)|P_{tt}(x_1,y-H(x_2))|\\
 &\ \ +H''(x_2)P_t(x_1,y-H(x_2))\geq0
 \,\ \ \text{in}\ \ I\times\R^2.
 \end{aligned}
\end{equation*}
Applying Kato's inequality \cite[Theorem~1.1]{BP2004},
we then obtain that
\begin{equation}\label{eq-epi-upper-cap}
 -\Delta
 P_H\geq\mathbf{1}_{\{P(x_1,y-H(x_2))\leq1\}}f(P(x_1,y-H(x_2)))
 =f(P_H)\,\ \ \text{weakly
 in}\ \ I\times\R^2.
\end{equation}
It follows from \eqref{eq-strip-front} that $0<W<B$ and $W\equiv0$ on
$\partial I\times\R$. Thus
$W(x_1,s_{-1})\leq\norm{W_{x_1}}_{L^\infty}\dist(x_1,\partial I)$. Since
$\max_{\overline I}B<1$, $\sup_{\R}\omega=1/(KR^2)$ and $\inf_{\overline
I\times\R}\Psi>0$, we obtain from \eqref{eq-epi-front-barriers} that for any
sufficiently large $R>0$,
\begin{equation}\label{eq-epi-side-negativity}
 \begin{split}
 &0<V_1\leq\max_{\overline
 I}B+\omega\norm{\Psi}_{L^\infty}<1\,\ \ \text{in}\ \ I\times\R^2,\\
& V_{-1}\leq2\omega^2\norm{W_{x_1}}_{L^\infty}-\omega\inf_{\overline
I\times\R}\Psi<0\,\ \ \text{in}\ \ \overline
I\times\R^2\cap\{\dist(x_1,\partial I)\leq2\omega(x_2)^2\}.
 \end{split}
\end{equation}
At $y=H(x_2)/2$, \eqref{eq-epi-modulated-fronts} and \eqref{eq-Hx2} imply that
$s_1\geq H(x_2)/8$. Since $P_H(x',H(x_2)/2)\leq P(x_1,-H(x_2)/2)$, we obtain
from \eqref{eq-strip-tails}, \eqref{eq-epi-cap-bounds} and
\eqref{eq-epi-front-barriers} that
\begin{equation*}
 (V_1-P_H)(x',H(x_2)/2)\geq\omega(x_2)\inf_{\overline I\times\R}\Psi-C(e^{-\mu
 H(x_2)/8}+e^{-kH(x_2)/2}).
\end{equation*}
Since $s_{-1}\leq-3H(x_2)/8$ for any $y\leq-H(x_2)$,
we deduce from \eqref{eq-strip-tails}
and \eqref{eq-epi-front-barriers} that
\[
 V_{-1}(x',y)
 \leq Ce^{-3\mu H(x_2)/8}
 -\omega(x_2)\inf_{\overline I\times\R}\Psi
 \ \ \text{for any}\ \ y\leq-H(x_2).
\]
Taking $c:=\min\{\mu/8,k/2,1/4\}$, we then obtain from
\eqref{eq-epi-front-parameters} and \eqref{eq-Hx2} that
\begin{equation}\label{eq-epi-gluing-error}
 \sup_{x_2\in\R}\frac{e^{-cH(x_2)}}{\omega(x_2)}
 =KR^2e^{-cH_0}\sup_{r\geq1}r^3e^{-4cr}\to0\ \ \text{as}\ \ H_0\to+\infty.
\end{equation}
This then implies that for any $x'\in\overline I\times\R$ and sufficiently
large $H_0>0$,
\begin{equation}\label{eq-epi-gluing-inequalities}
 V_1(x',H(x_2)/2)>P_H(x',H(x_2)/2),\ \ 
 V_{-1}(x',y)<0\ \ \text{for any}\ \ y\leq-H(x_2).
\end{equation}

Let $\Omega_3$ be the epigraph of
\begin{equation}\label{eq-epi-boundary-curve}
 \varphi(x_1,x_2):=H(x_2)+e^{H(x_2)}\frac{x_1^2-L^2}{1+x_1^2},
\end{equation}
where $H$ is given by \eqref{eq-Hx2}. To verify uniform smoothness, write
$p(s):=(s^2-L^2)/(1+s^2)$ and $\Phi(x',y):=e^{-H(x_2)}(y-H(x_2))-p(x_1)$. All
positive-order derivatives of $\Phi$ are bounded on $\{|\Phi|\leq1\}$. On
$\partial\Omega_3$ we have
\begin{equation*}
 \nabla\Phi(x',y)=\bigl(-p'(x_1),
 -H'(x_2)(p(x_1)+e^{-H(x_2)}),e^{-H(x_2)}\bigr).
\end{equation*}
Since $\inf_{\R}(|p|+|p'|)>0$ and $|H'(x_2)|\to4/R$ as $|x_2|\to\infty$, we
have $\inf_{\partial\Omega_3}|\nabla\Phi|>0$. We then obtain 
from the implicit function theorem
that $\Omega_3$ is uniformly smooth.
It follows from \eqref{eq-epi-boundary-curve} that
\[
 \varphi(x_1,x_2)\geq H(x_2)
 \ \ \text{for any}\ \ |x_1|\geq L,
\]
and hence
\[
 \Omega_3\cap\{y\leq H(x_2)\}\subset I\times\R^2.
\]

We may therefore define on $\Omega_3$,
\begin{equation}\label{eq-epi-final-barriers}
 \underline v(x',y):=\begin{cases}
 (V_{-1})_+,&x_1\in I,\\
 0,&x_1\notin I,
 \end{cases}
 \ \ 
 \overline v(x',y):=\begin{cases}
 \min\{V_1,P_H\},&y\leq H(x_2)/2,\\
 P_H,&H(x_2)/2<y<H(x_2),\\
 1,&y\geq H(x_2).
 \end{cases}
\end{equation}
Using the continuity of $V_1-P_H$ in $I\times\R^2$,
we then deduce from \eqref{eq-epi-gluing-inequalities}
that $\overline v=P_H$ in an open subset of $\Omega_3$
containing $\Omega_3\cap\{y=H(x_2)/2\}$. Since $P(x_1,0)>1$ on $\overline
I$ by \eqref{eq-epi-cap-bounds}, continuity of $P$ implies that $P_H\equiv1$ near
$y=H(x_2)$ in $I\times\R^2$ and hence $\overline v\equiv1$ near that surface
in $\Omega_3$.
For any $|x_1|<L$, \eqref{eq-epi-side-negativity} yields that $(V_{-1})_+=0$ if
$L-|x_1|\leq2\omega(x_2)^2$, and \eqref{eq-epi-gluing-inequalities} implies that
$(V_{-1})_+=0$ if $y\leq-H(x_2)$. For any sufficiently large $H_0>0$ and 
$(x',y)\in\partial\Omega_3$ satisfying
$|x_1|<L$ and $y>-H(x_2)$, we deduce from
\eqref{eq-epi-gluing-error} and \eqref{eq-epi-boundary-curve} that
\begin{equation*}
 L-|x_1|
 =\frac{1+x_1^2}{L+|x_1|}(H(x_2)-y)e^{-H(x_2)}
 \leq\frac{2(1+L^2)}L H(x_2)e^{-H(x_2)}
 \leq2\omega(x_2)^2.
\end{equation*}
 It follows from \eqref{eq-epi-side-negativity}
and \eqref{eq-epi-final-barriers} that
$\underline v\equiv0$ on $\partial\Omega_3$. We then
obtain from \eqref{eq-epi-strict-barriers}, \eqref{eq-epi-upper-cap},
\eqref{eq-epi-side-negativity} and \eqref{eq-epi-final-barriers} that for any
sufficiently large $H_0>0$,
\begin{equation}\label{eq-epi-subsuper}
 \begin{aligned}
&-\Delta\underline v\leq f(\underline v)\,\ \ \text{weakly
in}\ \ \Omega_3,\ \  0\leq\underline
v\leq1\,\ \ \text{in}\ \ \Omega_3,\ \ \underline
v\equiv0\,\ \ \text{on}\ \ \partial\Omega_3\\
&-\Delta\overline v\geq f(\overline v)\,\ \ \text{weakly
in}\ \ \Omega_3,\ \ 0<\overline
v\leq1\,\ \ \text{in}\ \ \Omega_3,\ \ \overline
v\geq0\,\ \ \text{on}\ \ \partial\Omega_3.
 \end{aligned}
\end{equation}

We finally construct two distinct positive bounded
solutions of \eqref{mainequation} by monotone iteration.
Fix $H_0$ so large that \eqref{eq-epi-gluing-inequalities} and
\eqref{eq-epi-subsuper} hold and $H_0 e^{-H_0}<L^2$. Choose a sufficiently
large ball $B_{R_0}(x_0)\Subset\Omega_3\cap\{y>H(x_2)\}$ and take
$w_0(x):=\tfrac12(1-|x-x_0|^2/R_0^2)_+^2$ satisfying
\eqref{eq-general-ball-subsolution}. We deduce from \eqref{eq-epi-final-barriers} that
$w_0\leq\overline v$ in $\Omega_3$. Let $u,U$ be the limits in
\eqref{eq-extremal-limits} for $T$ on $\Omega_3$. Applying
\eqref{eq-epigraph-ordered-iteration} to $(w_0,\overline v)$ and $(\underline
v,1)$ using \eqref{eq-epi-subsuper}, we then obtain that
\begin{equation*}
 w_0\leq u\leq\overline v,\ \ \underline v\leq U\leq1,\ \ u\leq
 U\,\ \ \text{in}\ \ \Omega_3.
\end{equation*}
Since $\varphi(0,x_2)=H(x_2)-L^2e^{H(x_2)}<0$, we have $(0,x_2,0)\in\Omega_3$.
It follows from \eqref{eq-epi-front-parameters} and
\eqref{eq-epi-modulated-fronts} that $s_1(x_2,0)\to-\infty$,
$s_{-1}(x_2,0)\to+\infty$ and $\omega(x_2)\to0$ as $|x_2|\to\infty$. We then
obtain from \eqref{eq-strip-tails}, \eqref{eq-epi-front-barriers} and
\eqref{eq-epi-final-barriers} that
\begin{equation}\label{eq-epi-barrier-separation}
 \underline v(0,x_2,0)\to B(0)>0,\ \ 
 0<\overline v(0,x_2,0)\leq V_1(0,x_2,0)\to0\ \ \text{as}\ \ |x_2|\to\infty.
\end{equation}
It follows from \eqref{eq-epi-barrier-separation} that
$u\leq\overline v<\underline v\leq U$ at $(0,x_2,0)$
for any sufficiently large $|x_2|$. The strong maximum
principle implies that $0<u<U<1$ in $\Omega_3$. For each $d>3$, extend $u,U$
constantly in the remaining $d-3$ horizontal coordinates on
$\{(x',y)\in\R^{d-1}\times\R:(x_1,x_2,y)\in\Omega_3\}$.
{This completes the proof of
Theorem~{\ref{thm-main}}(ii).}
\end{proof}

{
\begin{proof}[Proof of Theorem~{\ref{thm-extension}}]
Take $f,I,B,W$ from Lemma~\ref{lem-transverse-states}(iii).
Applying the construction in the proof of Theorem~{\ref{thm-main}}(ii),
we then obtain a uniformly $C^{2,\gamma}$ epigraph
$\Omega\subset\R^d$ on which \eqref{mainequation} admits
two positive bounded solutions $0<u<U<1$
for each integer $d\geq3$.
This completes the proof of Theorem~{\ref{thm-extension}}.
\end{proof}
}

\end{document}